\documentclass[11pt, a4paper]{article}

\usepackage{amsthm, amssymb, amsmath, textcomp, mathrsfs}
\usepackage[latin1,ansinew]{inputenc}
\usepackage[pdftex]{graphics}
\usepackage{graphicx}
\usepackage{longtable}
\usepackage{scrpage2}
\usepackage{hyphenat}
\usepackage{paralist}
\usepackage{color}

\makeatletter
\renewcommand{\section}{%
  \@startsection {section}{1}{\z@}%
                 {-3.5ex plus -1ex minus -.2ex}%
                 {2.3ex plus.2ex}%
                 {\normalfont\Large\bfseries}}

\renewcommand{\subsection}{%
  \@startsection {subsection}{1}{\z@}%
                 {-3.5ex plus -1ex minus -.2ex}%
                 {2.3ex plus.2ex}%
                 {\normalfont\large\bfseries}}
                 
\renewcommand{\subsubsection}{%
  \@startsection {subsubsection}{1}{\z@}%
                 {-2.5ex plus -1ex minus -.2ex}%
                 {1.3ex plus.2ex}%
                 {\normalfont\bfseries}}
\makeatother

\setheadsepline{.4pt}
\clearscrheadfoot
\lehead{\parbox[b]{0.8cm}{\pagemark}\textsc{Chapter \thechapter}}
\rohead{\textsc{\rightmark}\parbox[b]{0.8cm}{\hfill \pagemark}}

	\newcommand*{\savesymbol}[1]{%
	  \expandafter\let\csname orig#1\expandafter\endcsname\csname#1\endcsname
	  \expandafter\let\csname #1\endcsname\relax
	}
	
	\newcommand*{\restoresymbol}[2]{%
	  \expandafter\global\expandafter\let\csname#1#2\expandafter\endcsname%
	    \csname#2\endcsname
	  \expandafter\global\expandafter\let\csname#2\expandafter\endcsname%
	    \csname orig#2\endcsname
	}
	
	\newcommand*{\renamerobustsymbol}[2]{%
	  \expandafter\let\expandafter\origrealcommand
	    \csname #2\space\endcsname
	  \expandafter\global\expandafter\let\csname#1#2\endcsname=\origrealcommand
	}
	
	\makeatletter
	\def\ifnotsavedsym@helper#1#2!{\expandafter\ifx\csname orig#2\endcsname\relax}
	\newcommand*{\ifnotsavedsym}[1]{%
	  \expandafter\ifnotsavedsym@helper\string#1!%
	}
	\makeatother
	
	\newif\ifloadpackages
	\loadpackagestrue
	
	\makeatletter
	\newcommand{\missingpkgs}{}
	\newcommand{\foundpkgs}{}
	\newcommand{\if@sty@file@exists@star}[3]{%
	  \ifloadpackages
	    \IfFileExists{#1.sty}{#2}{#3}%
	  \else
	    #3%
	  \fi
	}
	\newcommand{\if@sty@file@exists}[3]{%
	  \ifloadpackages
	    \IfFileExists{#1.sty}%
	                 {#2\@cons\foundpkgs{{#1}}}%
	                 {#3\completefalse\@cons\missingpkgs{{#1}}}%
	  \else
	    #3\completefalse\@cons\missingpkgs{{#1}}%
	  \fi
	}
	\newcommand{\IfStyFileExists}{%
	  \@ifstar{\if@sty@file@exists@star}{\if@sty@file@exists}%
	}
	\makeatother

	\allowdisplaybreaks
	
	\usepackage{stmaryrd}
	
	\newif\ifMTOOLS
	
	\IfStyFileExists{mathtools}
  {\MTOOLStrue
   \savesymbol{xleftrightarrow} \savesymbol{xLeftarrow}
   \savesymbol{xRightarrow} \savesymbol{xLeftrightarrow}
   \savesymbol{xrightharpoondown} \savesymbol{xrightharpoonup}
   \savesymbol{xleftharpoondown} \savesymbol{xleftharpoonup}
   \savesymbol{xleftrightharpoons} \savesymbol{xrightleftharpoons}
   \savesymbol{xhookleftarrow} \savesymbol{xhookrightarrow}
   \savesymbol{xmapsto} \savesymbol{underbracket}
   \savesymbol{overbracket} \savesymbol{lparen} \savesymbol{rparen}
   \savesymbol{dblcolon} \savesymbol{coloneqq} \savesymbol{Coloneqq}
   \savesymbol{coloneq} \savesymbol{Coloneq} \savesymbol{eqqcolon}
   \savesymbol{Eqqcolon} \savesymbol{eqcolon} \savesymbol{Eqcolon}
   \savesymbol{colonapprox} \savesymbol{Colonapprox}
   \savesymbol{colonsim} \savesymbol{Colonsim} \savesymbol{overbrace}
   \savesymbol{underbrace}

   \let\origAtBeginDocument=\AtBeginDocument
   \def\AtBeginDocument##1{##1}
   \usepackage[donotfixamsmathbugs]{mathtools}
   \let\AtBeginDocument=\origAtBeginDocument

   \restoresymbol{MTOOLS}{xleftrightarrow}
   \restoresymbol{MTOOLS}{xLeftarrow}
   \restoresymbol{MTOOLS}{xRightarrow}
   \restoresymbol{MTOOLS}{xLeftrightarrow}
   \restoresymbol{MTOOLS}{xrightharpoondown}
   \restoresymbol{MTOOLS}{xrightharpoonup}
   \restoresymbol{MTOOLS}{xleftharpoondown}
   \restoresymbol{MTOOLS}{xleftharpoonup}
   \restoresymbol{MTOOLS}{xleftrightharpoons}
   \restoresymbol{MTOOLS}{xrightleftharpoons}
   \restoresymbol{MTOOLS}{xhookleftarrow}
   \restoresymbol{MTOOLS}{xhookrightarrow}
   \restoresymbol{MTOOLS}{xmapsto}
   \restoresymbol{MTOOLS}{underbracket}
   \restoresymbol{MTOOLS}{overbracket} \restoresymbol{MTOOLS}{lparen}
   \restoresymbol{MTOOLS}{rparen} \restoresymbol{MTOOLS}{dblcolon}
   \restoresymbol{MTOOLS}{coloneqq} \restoresymbol{MTOOLS}{Coloneqq}
   \restoresymbol{MTOOLS}{coloneq} \restoresymbol{MTOOLS}{Coloneq}
   \restoresymbol{MTOOLS}{eqqcolon} \restoresymbol{MTOOLS}{Eqqcolon}
   \restoresymbol{MTOOLS}{eqcolon} \restoresymbol{MTOOLS}{Eqcolon}
   \restoresymbol{MTOOLS}{colonapprox}
   \restoresymbol{MTOOLS}{Colonapprox}
   \restoresymbol{MTOOLS}{colonsim} \restoresymbol{MTOOLS}{Colonsim}
   \restoresymbol{MTOOLS}{overbrace} \restoresymbol{MTOOLS}{underbrace}

  }
  {}
  \newif\ifXPFEIL
	
	\IfStyFileExists{extpfeil}
	  {\XPFEILtrue
	   \let\origRequirePackage=\RequirePackage
	   \renewcommand*{\RequirePackage}[2][]{}
	   \savesymbol{xlongequal}
	   \savesymbol{xmapsto}
	   \savesymbol{xhookrightarrow}
	   \usepackage{extpfeil}
	   \restoresymbol{XPFEIL}{xlongequal}
	   \restoresymbol{XPFEIL}{xmapsto}	   
	   \let\RequirePackage=\origRequirePackage
	  }
	  {}

	\ifdefined\excludeHyperref
	\else
	\usepackage{hyperref}
	\hyperbaseurl{.}	
	\fi

\newtheorem*{theorem*}{Theorem}
\ifdefined\thmBookNumbering
	\newtheorem{theorem}{Theorem}[chapter]
	\newtheorem{definition}[theorem]{Definition}
	\newtheorem{remark}[theorem]{Remark}
	\newtheorem{corollary}[theorem]{Corollary}
	\newtheorem{proposition}[theorem]{Proposition}
	\newtheorem{lemma}[theorem]{Lemma}
	\newtheorem{example}[theorem]{Example}
	\newtheorem*{definition*}{Definition}
	\newtheorem*{maintheorem}{Main Theorem}
	\newtheorem*{unnumberedtheorem}{Theorem}
	\newtheorem*{remark*}{Remark}
	\newtheorem*{unnumberedcorollary}{Corollary}
\else
  \newtheorem{definition}{Definition}
	\newtheorem*{definition*}{Definition}
	\newtheorem{theorem}{Theorem}
	\newtheorem*{maintheorem}{Main Theorem}

	\newtheorem*{remark*}{Remark}

	\newtheorem{proposition}{Proposition}
	\newtheorem{lemma}{Lemma}
	
\fi

	\newcommand{\pr}{\operatorname{pr}}	
	\newcommand{\To}{\longrightarrow}
	
	\newcommand{\rank}{\operatorname{rank}}	
	
	\newcommand{\im}{\operatorname{im}}
	\newcommand{\Hol}{\operatorname{Hol}}

	\newcommand{\R}{\mathbb{R}}

	\newcommand{\Z}{\mathbb{Z}}

	\newcommand{\T}{\mathbb{T}}	
	\newcommand{\V}{\mathbb{V}}	
	\newcommand{\IS}{\mathbb{S}}
	\newcommand{\IL}{\mathbb{L}}
	
	\newcommand{\M}{\mathcal{M}}

	\newcommand{\cL}{\mathcal{L}}

	\newcommand{\Ric}{\operatorname{Ric}}

\begin{document}

	\title{\textsc{On Lorentzian manifolds with highest first Betti number}}
	  
  \author{\textsc{Daniel Schliebner}\thanks{The author is funded by the Berlin Mathematical School (BMS).}\\		
		\small \itshape Humboldt-Universit\"at zu Berlin, Institut f\"ur Mathematik, \\
		\small \itshape Rudower Chaussee 25, D--12489 Berlin, Germany.\\
		\small \itshape E-mail: schliebn@mathematik.hu-berlin.de, telephone: +49 (0)30 2093 1809.
	}

  \date{}

	\maketitle
	
  \begin{abstract}
  \noindent We consider Lorentzian manifolds with parallel light-like vector field $V$.
	Being parallel and light-like, the orthogonal complement of $V$ induces a codimension
	one foliation. Assuming compactness of the leaves and 
	non-negative Ricci curvature on the leaves it is known
	that the first Betti number is bounded by the dimension of the manifold or the
	leaves if the manifold is compact or non-compact, respectively.
	We prove in the case of the maximality of the first Betti number that every such Lorentzian manifold 
	is -- up to finite cover -- diffeomorphic to the torus (in the compact case)
	or the product of the real line with a torus (in the non-compact case) and has
	very degenerate curvature, i.e. the curvature tensor induced on the leaves is light-like.	   
	
  \medskip
  
  \end{abstract}	

	

	\section{Introduction}

%
%

We consider $(n + 2)$-dimensional Lorentzian manifolds\footnote{In this paper, all manifolds are assumed to be smooth, connected and without boundary.} $(\M, g)$ with parallel light-like vector field $V$, i.e.\ with
$g(V,V) = 0$ and $\nabla^gV = 0$ for the Levi-Civita connection $\nabla^g$ of $g$. Since
the vector field $V$ is parallel and light-like we obtain a parallel line bundle $\V := \R V$, 
while the orthogonal complement distribution
$$
	\V^\bot := \{ X \in \Gamma(T\M) \ | \ g(X,V) = 0 \}
$$
defines a parallel sub-distribution of the tangent bundle of codimension one and $\V \subset \V^\bot$.
Being parallel, $\V^\bot$ induces a codimension one foliation on $\M$ (which we will, in abuse of notation, also denote by $\V^\bot$) into light-like hypersurfaces 
which are the leaves of the foliation, each of which we usually will denote with $L$, i.e.\ $TL = \V^\bot$. The holonomy\footnote{By
	the \textit{holonomy} of $(\M, g)$ we mean the group
		$\Hol_x(\M, g) := \{ \mathcal P_\gamma \ | \ \gamma \text{ loop in } x \} \subset {\rm O}(T_x\M, g_x)$
	of parallel displacements along loops closed in $x \in \M$.
} 
of each such Lorentzian manifold is non-irreducible and contained in the stabilizer ${\rm O}(n) \ltimes \R^n$ of $V$ in ${\rm O}(1, n + 1)$.
If, moreover, the projection of the holonomy group onto $\R^n$ is surjective it is even indecomposable. An overview over this topic can be found
for example in \cite{galaev-leistner-esi,12}. A general interest in this area is to find topological consequences on the manifold $\M$ if it admits such a Lorentzian metric.
In light of this motivation, K. Lärz studied in his dissertation \cite{24}, among other aspect of global holonomy theory, topological questions on Lorentzian manifolds which
admit a parallel light-like vector field. One result he obtained was the boundedness of the first Betti number\footnote{We define the first Betti number of any manifold $\M$ to be the rank of $H^1(\M, \R)$.} $b_1(\M)$ by the dimension of $\M$ under additional assumptions on the
structure of the foliation and the Ricci curvature, c.f.\ \cite[Theorem 2.82]{24}, see also Lemma \ref{PropDimEstimate} in 
Section \ref{SecMain}. To be precise, he proved that every Lorentzian manifold with parallel light-like vector 
field\footnote{Essentially his assumption on the light-like vector field $V$ was weaker. Namely he just assumed $V$ to be \textit{recurrent}, 
i.e.\ such that $\nabla^gV = \omega \otimes V$, with $\omega \in \Omega^1(\M)$ and $\ker\omega = \V^\bot$.} such that the leaves of $\V^\bot$ are compact\footnote{Note that by
\cite{conlon1974transversally}, all leaves are diffeomorphic and either dense or closed. Indeed, a time-orientable indecomposable, non-irreducible Lorentzian manifold
is obviously transversally parallelizable, see e.g.\ \cite[Lemma 2.47]{24}.} and the Ricci curvature of 
$g$ is non-negative on $\V^\bot \times \V^\bot$,
fulfills $\varepsilon \leq b_1(\M) \leq \dim(\M) - 1 + \varepsilon$, where $\varepsilon = 1$  if $\M$ is compact, and $\varepsilon = 0$ if $\M$ is non-compact. \par
	A natural question occurring in this context is, if one can describe the cases in which the upper bound for the first Betti number is actually reached. This is basically motivated by the
classical Bochner result by which any compact, oriented Riemannian manifold $\mathcal N$ with non-negative Ricci curvature has $b_1(\mathcal N) \leq \dim \mathcal N$ and
$b_1(\mathcal N) = \dim \mathcal N$ if and only if it is isometric to the flat torus \cite[Ch. 7, Corollary 19]{petersen2006riemannian}. \par
	The main intention of this paper is to give a full answer to this question in the Lorentzian case under the assumptions that the \textit{leaves of $\V^\bot$ are compact}
and the \textit{Ricci curvature of $g$ is non-negative on $\V^\bot \times \V^\bot$}. Beside the ideas already used in \cite{24} our main contribution is to prove, along a leaf of $\V^\bot$, the existence of a 
\textit{horizontal} and \textit{integrable} realization of the screen bundle $\Sigma = \V^\bot/\V$, which is a subbundle $\IS$ of $T\M$ isomorphic to $\Sigma$. This
is done by exploiting the maximality assumption on the Betti number (Proposition~\ref{PropMain}) and it relates the Ricci curvature of $g$ along a leaf
with the Ricci curvature of a Riemannian metric (Lemma~\ref{LemHC}). We point out that horizontal and integrable screen distributions are a crucial tool in the study of Lorentzian manifolds with
indecomposable holonomy, see e.g.\ \cite{leistner05c, 24, leistnerschl2013complete}, since they provide a link between Riemannian and Lorentzian geometry.
	As a result we find out that under the assumptions above, the consequence for the Lorentzian metric $g$ is not to be flat
but still having very degenerate curvature. Indeed, the obtained metrics have \textit{light-like hypersurface curvature}. These metrics are in some sense 
generalizations of pp-waves and were in light of this motivation already studied by T. Leistner in \cite{leistner05c}.
Topologically, the manifolds turn out to be diffeomorphic (homeomorphic in dimension four) to a finite cover of the torus (in the compact case) or the product of the real line 
with the torus (in the non-compact case).

\begin{definition*}
	A Lorentzian manifold $(\M, g)$ with a global non-trivial parallel light-like vector field $V \in \Gamma(T\M)$, i.e.\
	$V \neq 0$, $g(V,V) = 0$ and $\nabla^gV = 0$, has {\rm{light-like hypersurface curvature}}, iff the curvature $R$ satisfies $R(X,Y)W \in \Gamma(\V)$ for all $X,Y,W \in \Gamma(\V^\bot)$.
\end{definition*}

\begingroup
\renewcommand\thetheorem{\Alph{theorem}}

\begin{maintheorem}
	\label{Thm-Main}
	Let $(\M, g)$ be an orientable $(n + 2)$-dimensional Lorentzian manifold with parallel light-like vector field $V$. Assume that the leaves of the codimension one
	foliation induced by the distribution $\V^\bot$ are compact and $\Ric|_{\V^\bot \times \V^\bot} \geq 0$. Then:	
	\begin{itemize}
		\item[(i)] if $\M$ is compact, then $b_1(\M) \leq n + 2$ and $b_1(\M) = n + 2$ if and only if $\M$ is -- up to finite cover -- diffeomorphic (homeomorphic\footnote{In dimension four, $\M$ is only 
		known to be \textit{homeomorphic} to the torus and to our best knowledge it seems to be an open problem in geometric topology 
		if $\M$ must be also diffeomorphic to the torus.} if $\dim\M = 4$) to the torus and $g$ has light-like hypersurface curvature;
		\item[(ii)] if $\M$ is non-compact, then $b_1(\M) \leq n + 1$ and $b_1(\M) = n + 1$ if and only if $\M$ is isometric to $\R \times \T^{n + 1}$ and $g$ has light-like hypersurface curvature.
	\end{itemize}
	In both cases, the leaves of $\V^\bot$ are all diffeomorphic to the torus $\T^{n + 1}$.
\end{maintheorem}


\noindent Concerning non-orientability of $\M$ we give the following

\begin{remark*}
	If $\M$ is non-orientable, the topological statements asserted in the Main Theorem inherit to the 2-fold orientation covering $\widehat{\M}$ for $\M$. In particular
	we can still conclude $b_1(\M) \leq n + 2$ since $b_1(\widehat{\M}) = b_1(\M) + b_{n - 1}(\M)$, c.f.\ \cite{brasher1969homology}. Moreover, since
	a covering is a local isometry, $(\M, g)$ has also light-like hypersurface curvature.
\end{remark*}

\noindent\textbf{Acknowledgements.} We would like to thank Helga Baum for motivating us to investigate the problem considered in this paper.

\endgroup

\section{Preliminaries}

Throughout this section, let $(\M, g)$ be a Lorentzian manifold with parallel light-like vector field. The techniques used within the present paper involve cohomology theory for foliations. 
Hence, we will present the required definitions and facts used in this paper but with just giving references to the statements in the literature. 
Indeed, the book \cite{tondeur1997geometry} is a very good reference for the most of the theory used here.

\subsection{Screen distributions and screen bundles}

By the parallel light-like vector field $V$ on $(\M, g)$, we obtain a filtration
$$
	\V \subset \V^\bot \subset T\M
$$
of the tangent bundle and a codimension two vector bundle
$$
	\Sigma := \V^\bot/\V \longrightarrow \M
$$
which is called the \textit{screen bundle} of $(\M, g)$. It comes with a naturally given connection $\nabla^\Sigma$, induced by the Levi-Civita connection $\nabla^g$
of $g$: $\nabla^\Sigma_X\varphi := [\nabla^g_XY]$ for some $Y \in \Gamma(\V^\bot)$ s.t.\ $[Y] = \varphi$, where $[ \cdot ] : \V^\bot \longrightarrow \Sigma$ is the projection.
Taking into account a non-canonical splitting $s : \Sigma \rightarrow \V^\bot$ of the short exact sequence 
$0 \longrightarrow \V \longrightarrow \V^\bot \longrightarrow \Sigma \longrightarrow 0$ one obtains a codimension two distribution $\IS := s(\Sigma) \subset T\M$,
called a \textit{screen distribution} or \textit{realization of the screen bundle} $\Sigma$. 
Each such realization of the screen bundle is in 1-to-1 correspondence with a \textit{screen vector field} $Z \in \Gamma(T\M)$ with
$g(Z,Z) = 0$ and $g(V,Z) = 1$.

\subsection{Riemannian foliations and Riemannian flows}
\label{subs:folandflow}

	Let us begin with mentioning that we do not need to distinguish between the leaves topologically. Since $\V^\bot$ is defined by the closed one form $\sigma := g(V,\cdot)$,
all leaves are diffeomorphic \cite[Corollary 3.31]{tondeur1997geometry} if $\M$ is closed.
Moreover, all leaves have trivial leaf holonomy (for a definition we refer to \cite[Section 3]{tondeur1997geometry}), c.f.\ \cite{reeb1950courbure} (see \cite[Theorem 3.29]{tondeur1997geometry} for an English version). Hence, if $\M$ is non-compact, then the assumption on all the leaves being compact implies that they are all diffeomorphic, \cite[Ch. 2, Corollary 8.6]{sharpe1997differential}.
In particular, in case \textit{(i)} of the Main Theorem we have that if one leaf $L$ of $\V^\bot$ is compact then so are all leaves. \par
	Given a foliated manifold $(\mathcal N, \mathcal F)$ and a Riemannian metric $h$ on $\cal N$, then $h$ is said to be \textit{bundle-like} for $(\mathcal N, \mathcal F)$,
iff $(\mathcal L_{X}h)(Z_1, Z_2) = 0$ for all $X \in \Gamma(T\mathcal F)$ and $Z_i \in \Gamma(T\mathcal F^{\bot_{h}})$, where $\mathcal L$ denotes the Lie-derivative.
In this case, $(\mathcal N, \mathcal F, h)$ is called a \textit{Riemannian foliation}. Given a foliated manifold $(\mathcal N, \mathcal F, h)$ with Riemannian metric 
$h$, consider the normal bundle $Q := T\mathcal N/T\mathcal F$ and the exact sequence $0 \longrightarrow T\mathcal F \longrightarrow T\mathcal N \stackrel{\pi}{\longrightarrow} Q \longrightarrow 0$. 
Then, for a splitting $s : Q \rightarrow T\mathcal F^{\bot_{h}}$, the metric $h$ induces a \textit{transversal metric} $h^T$ on $Q$ by $h^T := s^*h|_{T\mathcal F^{\bot_{h}}}$. We define a connection on $Q$ by
\begin{equation}
	\label{E-transverseLC-intro}
	\nabla^T_X\varphi := \begin{cases}
		\pi(\nabla^h_XY_\varphi), & X \in \Gamma(T\mathcal F^{\bot_h}), \\
		\pi([X,Y_\varphi]), & X \in \Gamma(T\mathcal F),
	\end{cases}
\end{equation}
for any $\varphi \in \Gamma(Q)$ and $s(\varphi) = Y_\varphi$. It is torsion-free \cite[Proposition 3.8]{tondeur1997geometry} and if $h$ is bundle-like, it is
metric \cite[Theorem 5.8]{tondeur1997geometry}. Moreover, if $h_1$ and $h_2$ are two bundle-like metrics w.r.t.\ $(\mathcal N, \mathcal F)$ such that
$h_1^T = h_2^T$, then $\nabla_1^T = \nabla_2^T$ \cite[Theorem 5.9]{tondeur1997geometry}. For bundle-like $h$, $\nabla^T$ is called the
\textit{transversal Levi-Civita connection} of the Riemannian foliation $(\mathcal N, \mathcal F, h)$. Given the transversal 
Levi-Civita connection $\nabla^T$ of a Riemannian foliation one obtains the corresponding curvature tensors $R^T$. 
Considering the bundle isomorphism $s(Q) \simeq Q$ one obtains the transversal Ricci curvature $\Ric^T$ defined through 
$\Ric^T(e_i, e_j) := \sum_{k = 1}^{\dim Q} R^T(e_i,e_k,e_k,e_j)$ for an $h^T$-orthonormal frame $\varphi_\ell \in \Gamma(Q)$ and $s(\varphi_\ell) = e_\ell$.
	
	Given a realization $\IS$ of the screen bundle, we obtain a naturally associated Riemannian metric $g^R$ corresponding to $\IS$ by 
\begin{equation}
	\label{E-RiemMetric}
	g^R(V,\cdot) := g(Z, \cdot), \ g^R(Z,\cdot) := g(V, \cdot), g^R(X,\cdot) := g(X, \cdot) \text{ for } X \in \Gamma(\IS)	
\end{equation}
and extension by linearity. This Riemannian metric is of certain interest since it relates aspects of Riemannian geometry naturally to the Lorentzian geometry. An important example
is the following: assume that the realization $\IS$ of $\Sigma$ is \textit{integrable} (i.e.\ $[\Gamma(\IS), \Gamma(\IS)] \subset \Gamma(\IS)$) and \textit{horizontal}
(i.e.\ $[\Gamma(\IS), \Gamma(\V)] \subset \Gamma(\IS)$) at least along $L$. Then $\nabla^g|_L$ (the induced connection on a leaf $L$) and $\nabla^R|_L$ (the Levi-Civita connection of the induced metric by $g^R$ on $L$)
only differ by a section $T\mathcal N \To \V$, $\nabla^R|_LV = 0$ and $\Ric^g|_{\V^\bot \times \V^\bot} = \Ric^{g^R}|_{\V^\bot \times \V^\bot}$. We will use this as a key step in the proof of our Proposition \ref{PropMain}.
Another important property of $g^R$ is the fact that it constitutes a \textit{Riemannian flow} on $L$. Let us explain this in more detail. 
A \textit{Riemannian flow} $(\mathcal N, \mathcal F, h)$ is a one-dimensional foliation $\mathcal F$ s.t.\ $h$ is bundle-like w.r.t.\ $\mathcal F$.
Now, since the parallel light-like vector field $V$ on $\M$ has no zero, it constitutes a one-dimensional foliation $\mathcal F$ on $\M$ and in particular on each leaf $L$ by its flow.
The important property for any $g^R$ as in \eqref{E-RiemMetric} associated to some realization $\IS$ of the screen bundle is now:

\begin{lemma}[{\cite[Lemma 2.48]{24}}]
	\label{LemRiemFlow}
	For each leaf $L$ of $\V^\bot$ and any realization $\IS$ of the screen bundle, $(L, \mathcal F, g^R)$ is a Riemannian flow, where $\mathcal F$ is given by the flow
	of $V$ restricted to $L$.
\end{lemma}

For a Riemannian flow $(\mathcal N, \mathcal F, h)$ defined by a non-singular vector field $V \in \Gamma(T\mathcal N)$ such that $h(V,V) = 1$, 
we define the \textit{mean curvature one form} $\kappa \in \Omega^1(\mathcal N)$ to be $\kappa := h(\nabla^h_VV, \cdot)$.
	
	Finally the following property of $g^R$ will be important:

\begin{lemma}
	\label{LemHC}
	Let $(\M, g)$ be a Lorentzian manifold with parallel light-like vector field $V$. If along $L$ there exists
	a horizontal and integrable realization $\IS$ of the screen bundle, then $\nabla^hV = 0$ and
	\begin{equation}
		\label{equ:diffR}
		[R^g(X,Y) - R^h(X,Y)]W \in \Gamma(\V) \text{ for all } X,Y,W \in \Gamma(\IS),
	\end{equation}
	where $h$ denotes the Riemannian metric on $L$ defined by $h(X,\cdot) := g(X, \cdot)$ for $X \in \Gamma(\IS)$,
	$h(V,V) = 1$ and extension by linearity. In particular we have that $\Ric^g|_{\V^\bot \times \V^\bot} = \Ric^h$.
\end{lemma}

\begin{proof}
	Since $h(V,V) = 1$,	we obtain $h(\nabla^h V, V) = 0$. Moreover, applying the horizontality and involutivity property of $\IS$, we see by the Koszul formula for $h$ that
	$$
		h(\nabla_{S_1}^h V, {S_2}) = g(\nabla_{S_1}^g V, {S_2}) = 0
	$$
	for all $S_1, S_2 \in \Gamma(\IS)$ since $\nabla^g V = 0$. This proves $\nabla^h V = 0$. Moreover, as 
	$$
		h(\nabla_{S_1}^h S_2, {S_3}) = g(\nabla_{S_1}^g S_2, {S_3})
	$$
	for all	$S_i \in \Gamma(\IS)$ we obtain $\nabla^g|_{\cL^\bot} {-} \nabla^h \in \Gamma(\IL)$ for all $X,Y \in \Gamma(T\cL^\bot)$
	which proves~\eqref{equ:diffR}.
\end{proof}

\subsection{Basic and twisted cohomology of foliations}

In this paragraph we deal with the term of basic cohomology. A comprehensive introduction can be found, e.g.\ in \cite[Chapter 4, Chapter 7]{tondeur1997geometry}.
Let $(\mathcal N, \mathcal F)$ be an arbitrary foliated manifold.

\begin{definition}
	A $k$-form $\alpha \in \Omega^k(\mathcal N)$ is called {\rm basic} iff $X \lrcorner \alpha = 0$ and $\mathcal L_X\alpha = 0$ for all $X \in \Gamma(T\mathcal F)$.
	The set of basic $k$-forms is denoted by $\Omega_B^k(\mathcal F)$.
\end{definition}

Since $\mathcal L_X d\alpha = d\mathcal L_X \alpha = 0$ and $X \lrcorner d\alpha = \mathcal L_X\alpha - d(X \lrcorner \alpha) = 0$ for all $X \in \Gamma(T\mathcal F)$ and
$\alpha \in \Omega_B^k(\mathcal F)$, we obtain a subcomplex $(\Omega_B^*(\mathcal F), d_B)$ of the de Rham-complex with $d_B := d|_{\Omega_B^*(\mathcal F)}$ and
hence a corresponding cohomology $H_B^*(\mathcal F)$, the\textit{ basic cohomology}. Given $d_B$ there is also a formal $L^2$-adjoint $\delta_B$
\cite[Theorem 7.10]{tondeur1997geometry} and hence a transversal Laplacian $\Delta_B = d_B\delta_B + \delta_Bd_B$ on $\Omega_B^*(\mathcal F)$. Naturally, we then say that 
a basic form $\alpha \in \Omega_B^k(\mathcal F)$ is \textit{basic-harmonic} iff $\Delta_B\alpha = 0$. For Riemannian flows we have the following
important result relating the mean curvature one form with the basic-harmonic forms:

\begin{theorem}[\cite{dominguez1998finiteness,mason2000application}]
	\label{Thm-basic-harmonic}
	Let $(\mathcal N, \mathcal{F}, g)$ be a Riemannian flow on a compact manifold $\mathcal N$. Then there exists a bundle-like metric $\widehat{g}$ on $\mathcal N$ such	
	that $\kappa$ is basic-harmonic and $g^T = \widehat{g}^T$.
\end{theorem}

Not enough, we need another definition of cohomology for foliations. Namely, since for the basic cohomology, Poincar{\'e} duality not necessarily holds,
one defines another cohomology theory. The so called \textit{twisted cohomology} for a foliation $(\mathcal N, \mathcal F, h)$ with $\kappa \in \Omega_B^1(\mathcal N)$ is defined by the subcomplex
$(\Omega_B^*(\mathcal F), d_\kappa := d_B - \kappa \wedge \cdot)$ of the de Rham-complex.\footnote{Here, $\kappa$ is in general the mean curvature form of
the foliation \cite[(3.20)]{tondeur1997geometry}. However, we will use this only in the setting of Riemannian flows s.t.\ $\kappa$ is the mean curvature one form.} 
We will denote its cohomology groups by $H_\kappa^*(\mathcal F)$.
It behaves in a nice way as to obtain a type of Poincar{\'e} duality when comparing the basic cohomology and the twisted cohomology \cite[Theorem 7.54]{tondeur1997geometry}.

The orientability assumption in the Main Theorem has the following background.
Since most of the results about Hodge theory in basic cohomology require the foliation to be transversally oriented (i.e.\ there is an orientation of the
normal bundle), we need the following

\begin{lemma}
	\label{LemOrient}
	Let $(\M, g)$ be an oriented Lorentzian manifold with parallel light-like vector field. Then, for each leaf $L$ of $\V^\bot$ and any realization $\IS$ of the screen bundle, 
	the Riemannian flow $(L, \mathcal F, g^R)$ is transversally orientable. Moreover, $L$ is orientable.
\end{lemma}

\begin{proof}
	The argument can be found in \cite[p. 78]{24}. Namely, one can prove that $\Hol(\nabla^T) \subset {\rm SO}(\dim\IS)$ and hence $\Sigma$ is orientable. In particular,
	each leaf $L$ is orientable since any screen vector field $Z \in \Gamma(T\M)$ defines a unit normal vector field.
\end{proof}

\section{Main Results}
\label{SecMain}

Let $(\M, g)$ be an oriented Lorentzian manifold with parallel light-like vector field $V$ inducing an integrable codimension one distribution $\V^\bot \subset T\M$
and hence a codimension one foliation on $\M$. Recall from the previous section that, by fixing a realization $\IS$ of the screen bundle and denoting by 
$g^R$ the associated Riemannian metric on $\M$, 
we obtain a Riemannian flow $\mathcal F$ on $L$ by the flow of $V$ and $g^R$ restricted to $L$ (Lemma \ref{LemRiemFlow}). Throughout this section we will denote
this Riemannian flow by $(L, \mathcal F, g^R)$. 

\begin{lemma}
	\label{LemDirectSum}
	Let $L$ be a compact $(n + 1)$-dimensional leaf of $\V^\bot$. Then $H_{\rm dR}^1(L) = H_B^1(\mathcal F) \oplus H$ for 
	$H$ a subgroup of $H_B^{n}(\mathcal F) \in \{ 0, \R\}$.
\end{lemma}

\begin{proof}
	Since $(L, \mathcal F, g^R)$ is a Riemannian flow, by Theorem \ref{Thm-basic-harmonic}, there exists a bundle-like metric $\widehat{g}^R$ such that the
	mean-curvature one form $\kappa$ of $(L, \mathcal F, \widehat{g}^R)$ is basic-harmonic. By the Gysin sequence for $(L, \widehat{g}^R)$, 
	c.f.\ \cite[Theorem 3.2]{prieto2001gysin}, we obtain an exact sequence
	\begin{equation}
		\label{E-gysin}
		0 \longrightarrow H_B^1(\mathcal F) \longrightarrow H_{\rm dR}^1(L) \longrightarrow H_{\kappa}^0(\mathcal F) \longrightarrow H_B^2(\mathcal F) \longrightarrow \ldots
	\end{equation}
	and by taking into account that $H_{\kappa}^0(\mathcal F) \cong H_B^{n}(\mathcal F)$, c.f.\ \cite[Theorem 7.54]{tondeur1997geometry},
	this translates into the long exact sequence
	\begin{equation}
		\label{E-gysin-2}
		0 \longrightarrow H_B^1(\mathcal F) \longrightarrow H_{\rm dR}^1(L) \stackrel{\Phi}{\longrightarrow} H_B^{n}(\mathcal F) \longrightarrow H_B^2(\mathcal F) \longrightarrow \ldots
	\end{equation}
	and thus we obtain the short exact sequence
	\begin{equation}
		\label{E-short}
		0 \longrightarrow H_B^1(\mathcal F) \longrightarrow H_{\rm dR}^1(L) \longrightarrow H \longrightarrow 0
	\end{equation}
	for $H:= \im \Phi \subset H_B^{n}(\mathcal F)$. In particular this is a short exact sequence
	of vector spaces and hence splits as a direct sum. 
	Since $H_B^{n}(\mathcal F) \in \{ 0, \R\}$ by \cite[Corollary 7.57]{tondeur1997geometry} this completes the proof.
\end{proof}

By assuming non-negativity of the Ricci curvature on $TL \times TL$, we obtain the following estimation for the dimensions
of $H_{\rm dR}^1(L)$ and $H_B^1(\mathcal F)$. This is precisely the \cite[Proposition 2.81]{24} of K. L{\"a}rz. However, to be self-contained,
make the upcoming proofs more precise and to fix notation we will present its full proof here.

\begin{lemma}
	\label{PropDimEstimate}
	Let $L$ be a compact $(n + 1)$-dimensional leaf of $\V^\bot$ and $\Ric^g|_{TL \times TL} \geq 0$. Then 
	$b_1(L) \leq \dim H_B^1(\mathcal F) + 1 \leq \dim L$.
\end{lemma}

\begin{proof}
	Again, consider the Riemannian flow $(L, \mathcal F, g^R)$ and the bundle-like metric $\widehat{g}^R$ such that the
	mean-curvature one form $\kappa$ of $(L, \mathcal F, \widehat{g}^R)$ is basic-harmonic, constituted by Theorem \ref{Thm-basic-harmonic}.
	In particular, the induced metrics on $L/T\mathcal F =: \Sigma|_L$ coincide and hence so do the induced transversal connections.	
	Defining by $\widehat{\IS} := \ker\widehat{g}^R(V, \cdot)$ a realization of $\Sigma|_L$, the transversal Levi-Civita connection
	$\nabla^T : \Gamma(\widehat{\IS}) \longrightarrow \Gamma(T^*L \otimes \widehat{\IS})$ is given by
	\begin{equation}
		\label{E-transverseLC}
		\nabla^T_XY := \begin{cases}
			\pr_{\widehat{\IS}}(\widehat{\nabla}_XY), & X \in \Gamma(\widehat{\IS}), \\
			\pr_{\widehat{\IS}}([X,Y]), & X \in \Gamma(\V),
		\end{cases}
	\end{equation}
	for any $Y \in \Gamma(\widehat{\IS})$, where $\widehat{\nabla}$ is the Levi-Civita connection of $\widehat{g}^R$ (see Section \ref{subs:folandflow}). 
	Since $\Ric^g(V, \cdot) = 0$ and $g_{\widehat{\IS} \times \widehat{\IS}} = g^R|_{\widehat{\IS} \times \widehat{\IS}}$ we have that
	\begin{equation}
		\label{E-Ric}
		\Ric^T_{\widehat{\IS} \times \widehat{\IS}} = \Ric^g|_{\widehat{\IS} \times \widehat{\IS}} \geq 0.
	\end{equation}
	Since any class in $H_B^1(\mathcal F)$ can be represented by a $\Delta_B$-harmonic one-form $\alpha \in \Omega_B^1(L)$ 
	\cite[Theorem 7.51]{tondeur1997geometry}, we choose such a basic-harmonic $\alpha$ for each generator of $H_B^1(\mathcal F)$.
	To apply a Bochner argument we need an appropriate Weizenb{\"o}ck formula which is given by
	\cite[Proposition 6.7]{habib2010modified} and by integration reads
	$$
		0 = \int_L || \nabla^T \alpha ||^2 + \int_L \Ric^T(\alpha^\sharp, \alpha^\sharp),
	$$
	where $\alpha^\sharp$ is the dual of $\alpha$ w.r.t.\ $\widehat{g}^R$.
	Hence $\nabla_X^T \alpha = 0$ for all $X \in \Gamma(TL)$. This proves $\dim H_B^1(\mathcal F) \leq n = \dim\widehat{\IS}$
	and thus the second asserted inequality. For the first inequality we apply Lemma \ref{LemDirectSum};
	this completes the proof.
\end{proof}

We can now prove the key step for the proof of the main result. Namely, the -- by the previous lemma -- maximality of 
the dimension of the first Betti number of $L$ and $\Ric^g|_{TL \times TL} \geq 0$ already yield that
$L$ is the torus with the curvature of $\nabla^g|_{L}$ on $L$ being light-like.

\begin{proposition}
	\label{PropMain}
	Let $(\M^{n + 2}, g)$ be an oriented $(n + 2)$-dimensional Lorentzian manifold with parallel light-like vector field $V$ such that
	$\Ric^g|_{TL \times TL} \geq 0$ and the leaves $L$ of the codimension one foliation induced by $\V^\bot$ are compact.
	Then, $b_1(L) \leq n + 1$ and $b_1(L) = n + 1$ if and only if $\nabla^g|_{L}$ has light-like curvature and $L$ is diffeomorphic to the torus.
\end{proposition}

\begin{proof}
	The proof works as follows: Let $b_1(L) = n + 1$. We will define for every realization of the screen bundle another screen distribution
	\textit{along $L$} which is horizontal and integrable. The associated Riemannian metric $h$ to this screen distribution
	then has non-negative Ricci curvature by Lemma \ref{LemHC} and hence $(L, h)$ is the flat torus implying $\nabla^g|_{L}$ to have light-like curvature. \par	
	
		\textbf{Step 1.} Fix a realization of the screen bundle and hence a Riemannian flow $(L, \mathcal F, g^R)$ on $L$
	and denote by $\widehat{g}^R$ the bundle-like metric s.t.\ the mean-curvature one form $\kappa$ is basic-harmonic.
	W.l.o.g.\ we may additionally assume $\widehat{g}^R(V,V) = 1$. We define $\chi := \widehat{g}^R(V, \cdot)$. \par
	
		\textbf{Step 2.} Recall that $H_B^{n}(\mathcal F) \in \{ 0, \R\}$ (see Lemma \ref{LemDirectSum}) since $(L, \mathcal F, g^R)$ is transversally orientable by Lemma \ref{LemOrient}.		
		Assume that $H_B^{n}(\mathcal F) = 0$.
	Then we would have by the sequence \eqref{E-gysin-2} that $H_B^1(\mathcal F) \cong  H_{\rm dR}^1(L)$ which is impossible due to $b_1(L) = n + 1$
	and Lemma \ref{PropDimEstimate}. Hence $H_B^{n}(\mathcal F) = \R$ and \cite[Theorem A]{molino1985deux} 
	(see also \cite[Theorem 6.17]{tondeur1997geometry} for an English version) implies that
	$(L, \mathcal F, \widehat{g}^R)$ is an \textit{isometric} Riemannian flow, i.e.\ $V$ is a $\widehat{g}^R$-Killing field and therefore
	$\mathcal L_V \chi = 0$ and $\kappa = 0$. \par	
	
		\textbf{Step 3.} By Lemma \ref{PropDimEstimate} and $b_1(L) = n + 1$ we infer $\dim H_B^1(\mathcal F) = n$. For $(L, \mathcal F, \widehat{g}^R)$
	being an \textit{isometric} Riemannian flow, the Gysin long exact sequence \eqref{E-gysin-2} is given as
	\begin{equation}
		\label{E-gysin-3}
		0 \longrightarrow H_B^1(\mathcal F) \longrightarrow H_{\rm dR}^1(L) \stackrel{\Phi}{\longrightarrow} H_B^{0}(\mathcal F) \stackrel{\delta}{\longrightarrow} H_B^2(\mathcal F) \longrightarrow \ldots
	\end{equation}
	where $\Phi = (V \lrcorner \cdot)_*$ and $\delta = [d\chi \wedge \cdot]$, see \cite[Theorem 6.13]{tondeur1997geometry}.
	Note that since $H_{\kappa}^0(\mathcal F) \cong H_B^{n}(\mathcal F)$ by \cite[Theorem 7.54]{tondeur1997geometry} and
	$\kappa = 0$ we have $H_B^{0}(\mathcal F) = H_{\kappa}^0(\mathcal F) \cong H_B^{n}(\mathcal F) = \R$.
	This together with \eqref{E-gysin-3} implies the short exact sequence 
	\begin{equation}
		\label{E-short-3}
		0 \longrightarrow H_B^1(\mathcal F) \longrightarrow H_{\rm dR}^1(L) \longrightarrow \im \Phi \longrightarrow 0
	\end{equation}
	with $\im \Phi = \R$ since $b_1(L) = n + \dim \im \Phi$.	
	Hence, by the exactness of \eqref{E-gysin-3}, $\ker\delta = \R$ and so $0 = \delta([1]) = [1 \cdot d\chi]$. Hence $[d\chi]$ vanishes in $H_B^2(\mathcal F)$.
	In this case, $d\chi = d_B\alpha$ for some $\alpha \in \Omega_B^1(L)$. Then, $\omega := \chi - \alpha \in \Omega^1(L)$ is closed w.r.t.\ $d$ and $\omega(V) = 1$. We define
	\begin{equation}
		\label{E-new-screen-2}
		\IS := \ker\omega.
	\end{equation}
	Obviously, $\IS$ is integrable and it is horizontal since $\mathcal L_V\omega = \mathcal L_V \chi - \mathcal L_V\alpha = 0$ as $\alpha$ is a basic 1-form
	and $\mathcal L_V \chi = 0$ by Step 2.\par
	
	\textbf{Step 4.} Define a Riemannian metric on $L$ associated to $g$ and $\IS$ by 
	$$
		h(X,Y) := 
		\begin{cases}
			1, & X = Y = V \\
			g(X,Y), & X,Y \in \Gamma(\IS) \\
			0, & X \in \Gamma(\V) \text{ and } Y \in \Gamma(\IS) \text{ or } Y \in \Gamma(\V) \text{ and } X \in \Gamma(\IS)
		\end{cases}
	$$
	and linear extension. Since $\IS$ is horizontal and integrable, Lemma \ref{LemHC} yields
	$\Ric^h = \Ric^g|_{TL \times TL} \geq 0$. Therefore, $(L, h)$ turns into a compact, orientable Riemannian manifold
	with non-negative Ricci-curvature and is thus isometric to the flat torus \cite[Ch. 7, Corollary 19]{petersen2006riemannian}.
	In particular, equation \eqref{equ:diffR} in Lemma \ref{LemHC} implies that $\nabla^g|_{L}$ has light-like curvature.
\end{proof}

The last lemma of this section now relates the dimension of $H_{\rm dR}^1(\M)$ with the dimension of
$H_{\rm dR}^1(L)$ for a compact leaf of the codimension one foliation $\V^\bot$.

\begin{lemma}
	\label{PropCohomology}
	Let $(\M, g)$ be an $(n + 2)$-dimensional Lorentzian manifold with parallel light-like vector field $V$ such that
	the leaves $L$ of the codimension one foliation induced by $\V^\bot$ are compact. Then, 
	
	\begin{itemize}
		\item[(i)]  $b_1(\M) \leq b_1(L) + 1$.
		\item[(ii)]  $\pi_1(\M) \cong \Z \prescript{}{\varphi}{\ltimes} \pi_1(L)$ for some homomorphism $\varphi \in \operatorname{Hom}(\Z, \operatorname{Aut}(\pi_1(L)))$ and hence
		$H_1(\M, \Z) \cong \Z \oplus H_1(L, \Z)/K$, for the subgroup $K$ of $\pi_1(L)$ generated by the elements $\varphi(k)([\gamma]) \cdot [\gamma]^{-1}$, $k \in \Z$, $[\gamma] \in \pi_1(L)$.
	\end{itemize}
\end{lemma}

\begin{proof}
	Since $\V^\bot$ is defined as $\V^\bot = \ker \sigma$ for $\sigma := g(V, \cdot)$ and $d\sigma = 0$, all leaves have trivial leaf holonomy (see Section \ref{subs:folandflow}).	
	The compactness assumption for the leaves then implies, that $\M$ fibers over $S^1$ or $\R$ with each leaf being given as a fiber of the fibration \cite[Corollary 8.6]{sharpe1997differential}. Hence, if $\M$ is not compact,
	$\M$ fibers over $\R$ and thus $\M \simeq \R \times L$ for a fixed leaf. Otherwise, if $\M$ is compact, $\M$ fibers over $S^1$ and the long exact sequence of homotopy groups
	for fibrations \cite[Theorem 4.41]{hatcher2002algebraic} yields
	$$
		0 \longrightarrow \pi_1(L) \longrightarrow \pi_1(\M) \longrightarrow \Z \longrightarrow 0
	$$
	and hence $\pi_1(\M) \cong \Z \prescript{}{\varphi}{\ltimes} \pi_1(L)$ for some $\varphi \in \operatorname{Hom}(\Z, \operatorname{Aut}(\pi_1(L)))$ since $\Z$ is free. By the Hurewicz theorem
	we have $H_1(\M, \Z) \cong \pi_1(\M)_{\rm ab}$, where for any group $G$ we denote with $G_{\rm ab} := G/[G,G]$ its Abelization. 
	However, for any semi-direct product $G = A \prescript{}{\psi}{\ltimes} B$ for some $\psi \in \operatorname{Hom}(A, \operatorname{Aut}(B))$ we have that
	$$
		G_{\rm ab} = A_{\rm ab} \oplus B_{\rm ab}/H
	$$
	with $H$ denoting the subgroup of $B$ generated by elements $\psi(a)b\cdot b^{-1}$ with $a \in A$ and $b \in B$, c.f.\ \cite[Proposition 3.3]{gonccalves2009lower}.
	Therefore, we see that 
	$$
		H_1(\M, \Z) \cong \pi_1(\M)_{\rm ab} \cong \Z \oplus H_1(L, \Z)/K
	$$
	for the subgroup $K$ generated by the elements $\varphi(k)([\gamma]) \cdot [\gamma]^{-1}$, $k \in \Z$, $[\gamma] \in \pi_1(L)$. In particular,
	$b_1(\M) = 1 + \rank H_1(L, \Z) - \rank K \leq b_1(L) + 1$.
\end{proof}

\section{Proof of the Main Theorem}


Let us first proof part \textit{(i)} of the Main Theorem. \par 
	With $b_1(\M) = n + 2$ we infer $b_1(L) = n + 1$ by Lemma \ref{PropCohomology}\textit{(i)} taking into account
that $b_1(L) \leq \dim L = n + 1$ by Lemma \ref{PropDimEstimate}. Hence, all assumptions of Proposition \ref{PropMain} are satisfied. We obtain, that
the connection $\nabla^g|_L$ on each leaf $L$ induced by the Levi-Civita connection of $g$ has light-like curvature. By \cite[Proposition 6]{leistner05c} this is equivalent for $(\M, g)$
to have light-like hypersurface curvature. \par
Since $L = \T^{n + 1}$ and hence $\M$ fibers over $S^1$ with toric fibers, the long exact sequence of homotopy groups for fibrations implies that for $k > 1$ all homotopy groups $\pi_k(\M)$ vanish.
Therefore, $\M$ is a $K(\pi_1(\M), 1)$-space\footnote{With $K(G, 1)$-spaces we denote the Eilenberg–MacLane spaces with first fundamental group equal to $G$ and all other homotopy groups vanishing.},
while $\pi_1(\M) \cong \Z \prescript{}{\varphi}{\ltimes} \pi_1(L) \cong \Z \prescript{}{\varphi}{\ltimes} \Z^{n + 1}$ by Lemma \ref{PropCohomology}\textit{(ii)} and since
$\pi_1(L) = \Z^{n + 1}$. By assumption, $H_1(\M, \Z) = \Z^{n + 2} \oplus \operatorname{Tor}$ and hence
$$
	\Z^{n + 2} \oplus \operatorname{Tor} = H_1(\M, \Z) = \Z \oplus \Z^{n + 1}/K
$$
by Lemma \ref{PropCohomology}\textit{(ii)}. But this equation can only hold, if and only if $K$ is trivial. Namely, comparing the ranks of the left and
the right hand side, we observe
$$
	n + 2 = 1 + (n + 1) - \rank K \Longleftrightarrow \rank K = 0.
$$
As $K$ is a subgroup of $\pi_1(L) \cong \Z^{n + 1}$, $K$ is trivial. Therefore, $\M^{n + 2}$ is a $K(\Z^{n + 2}, 1)$-space and hence
homotopy-equivalent to the torus \cite[Theorem 1B.8]{hatcher2002algebraic}. But in the case of the torus, for $\dim\M \leq 3$ this is even equivalent for $\M$ itself\footnote{
For the 3-dimensional case, \textsc{Waldhausen} \cite{waldhausen1968irreducible} proved this for \textit{Haken manifolds} and thus in particular for 3-dimensional closed manifolds fibering over the circle. Note that
for 3-dimensional manifolds classifications up to diffeomorphism and homeomorphism coincide \cite{moise1952affine}.}
or some finite cover (if $\dim\M > 4$) to be diffeomorphic to the standard torus \cite[Page 236]{wall1999surgery}. For $\dim \M = 4$ we can only conclude that $\M$ is homeomorphic to $\T^4$ \cite[Chapter 11.5]{freedman2014topology}. \par
	It remains to prove part \textit{(ii)}. This is straightforward since $\M \simeq \R \times L$ for a fixed leaf $L$ by \cite[Corollary 8.6]{sharpe1997differential}
(see the proof of Lemma \ref{PropCohomology} for details) and hence $n + 1 = b_1(\M) = b_1(L)$ so	Proposition \ref{PropMain} applies.
\hfill $\Box$

	
		
	\small

\bibliographystyle{amsalpha}
\bibliography{Bibliography}

\end{document}